\documentclass[a4paper,10pt]{article}

\usepackage{amsmath,amsfonts,mathrsfs,amssymb,latexsym}
\usepackage[francais,english]{babel}
\usepackage[top=4cm, bottom=3cm, left=3cm, right=3cm]{geometry}

\newtheorem{theorem}{Theorem}

\newtheorem{e-proposition}[theorem]{Proposition}

\newtheorem{e-definition}[theorem]{Definition\rm}
\newtheorem{remark}[theorem]{\it Remark\/}

\newtheorem{theoreme}{Th\'eor\`eme}

\newcommand{\D}{{\mathrm{d}}}

\newcommand{\bJ}{\mathbf{J}}
\newcommand{\bH}{\mathbf{H}}

\begin{document}

\title{A one-dimensional Keller-Segel equation with a drift issued from the boundary}

\author{Vincent Calvez
\protect\footnote{Unit\'e de Math\'ematiques Pures et Appliqu\'ees, CNRS UMR 5669, \'Ecole Normale Sup\'erieure de Lyon, 46 all\'ee d'Italie, 69364 Lyon Cedex 07, France. E-mail: 
\texttt{vincent.calvez@umpa.ens-lyon.fr}
}\, ,
Nicolas Meunier
\protect\footnote{MAP5, Universit\'e
Paris Descartes, 45 rue des Saints-P\`eres, 75270 Paris cedex
06, France.
E-mail: \texttt{nicolas.meunier@parisdescartes.fr}
}
}

\maketitle

\begin{abstract}
We investigate in this note the dynamics of a one-dimensional Keller-Segel type model on the half-line. On the contrary to the classical configuration, the chemical production term is located on the boundary. We prove, under suitable assumptions, the following dichotomy which is reminiscent of the two-dimensional Keller-Segel system. Solutions are global if the mass is below the critical mass, they blow-up in finite time above the critical mass, and they converge to some equilibrium at the critical mass. Entropy techniques are presented which aim at providing quantitative convergence results for the subcritical case. This note is completed with a brief introduction to a more realistic model (still one-dimensional). 
\end{abstract}

\selectlanguage{francais} 
\begin{abstract}
Nous \'etudions dans cette note la dynamique d'un mod\`ele unidimensionnel de type Keller-Segel pos\'e sur une demi-droite. Dans le cas pr\'esent, la production du signal chimique est localis\'ee sur le bord, au lieu d'\^etre r\'epartie \`a l'int\'erieur du domaine comme dans le cas classique. On d\'emontre, sous des hypoth\`eses convenables, la dichotomie suivante qui rappelle le syst\`eme de Keller-Segel en dimension deux d'espace. Les solutions sont globales si la masse est sous-critique, elles explosent en temps fini si la masse d\'epasse la masse critique. Enfin, les solutions convergent vers un \'etat d'\'equilibre lorsque la masse est \'egale \`a la valeur critique. Des m\'ethodes d'entropie sont d\'evelopp\'ees, dans le but d'obtenir des r\'esultats de convergence quantitatifs. Cette note est enrichie d'une br\`eve introduction \`a un mod\`ele plus r\'ealiste (\`a nouveau unidimensionnel).
\end{abstract}

\selectlanguage{francais}
\section*{Version fran\c caise abr\'eg\'ee}

Dans cette note nous allons \'etudier le comportement math\'ematique en dimension un de l'\'equation aux d\'eriv\'ees partielles suivante :
\begin{equation}\label{eqfran}
\partial _t n(t,x) = \partial _{xx} n(t,x) +n(t,0) \partial _x n(t,x)\, , \quad t >0\, , \, x\in (0,+\infty)\, ,
\end{equation}
avec la condition initiale : $n(t=0,x) = n_0(x)\geq 0$.
Nous imposons au bord une condition de flux nul~:
$\partial _x n (t,0)+n(t,0)^2=0$, de sorte que la  masse est conserv\'ee au cours du temps (au moins formellement)~:
\begin{equation}\label{massconsfran}
\int_{x>0} n(t,x)\,\D x =\int_{x>0} n_0(x)\,\D x = M\, . 
\end{equation}

Ce mod\`ele a \'et\'e propos\'e dans \cite{Voituriez} pour d\'ecrire synth\'etiquement la polarisation des cellules de levure.
Une caract\'eristique int\'eressante de (\ref{eqfran}) r\'eside dans le fait que la solution peut devenir non born\'ee en temps fini. Dans cette note nous allons montrer l'alternative suivante~:\medskip

\begin{theoreme}[Existence globale {\em vs.} explosion] Supposons que $n_0(x)$ est continue sur $[0,+\infty)$ et que $n_0\in L^1_+((1+x)dx)$. 
Si $M\leq1$ alors la solution de \eqref{eqfran} est globale en temps. Au contraire si $M>1$, en supposant en outre que  $n_0$ est d\'ecroissante, alors la solution de  \eqref{eqfran}  explose en temps fini.
\end{theoreme}\medskip

Nous annon\c{c}ons \'egalement les r\'esultats suivants concernant le comportement asymptotique de la solution lorsque $M\leq 1$~:\medskip
\begin{theoreme}[Comportement asymptotique] Dans le cas critique $M = 1$, il existe une famille d'\'etats stationnaires pour \eqref{eqfran}  param\'etr\'ee par $\alpha>0$. La solution converge (au sens de l'entropie relative \eqref{eq:relative entropy}) vers l'\'equilibre tel que $\alpha^{-1} = \int_{x>0} x n_0(x)\, \D x$.\\
Dans le cas sous-critique $M<1$, la solution d\'ecro\^it vers z\'ero, et converge (au sens de l'entropie relative) vers un unique profil auto-similaire.  
\end{theoreme}\medskip

Enfin, nous nous int\'eressons \`a l'\'etude d'un mod\`ele plus r\'ealiste, qui prend en compte l'\'echange entre des particules libres, et des particules fix\'ees au bord qui cr\'eent le potentiel attractif (concentration $\mu(t)$)~:
\begin{equation*} 
\left\{\begin{array}{l} 
\partial _t n (t,x)=\partial _{xx} n (t,x) + \mu (t) \partial _x n (t,x) \, , \quad t >0\, , \, x\in (0,+\infty)\, , \\ 
 \mu'(t) = n(t,0)- \mu(t)\, ,
\end{array}\right.
\end{equation*}
avec la condition de flux au bord~: $\partial _x n (t,0)+\mu(t) n (t,0)=\mu'(t)$. 
\medskip

\begin{theoreme}
Avec les hypoth\`eses des th\'eor\`emes pr\'ec\'edents, et dans le cas sur-critique $M>1$, $\mu(t)$ converge vers $\overline{\mu} = M-1$ et la densit\'e $n(t,x)$ converge en entropie relative vers $h(x) = \overline{\mu} \exp\left(-\overline{\mu} x\right)$.
\end{theoreme}\medskip

\selectlanguage{english}
\section*{English Version}
In this note we shall study the mathematical behavior of the following one dimensional partial differential equation:
\begin{equation}\label{eq}
\partial _t n(t,x) = \partial _{xx} n(t,x) +n(t,0) \partial _x n(t,x)\, , \quad t >0\, , \, x\in (0,+\infty)\, ,
\end{equation}
together with the initial condition: $n(t=0,x) = n_0(x)\geq 0$.
We impose a zero-flux boundary condition for the density $n$,
\begin{equation}\label{cl}
\partial _x n (t,0)+n(t,0)^2=0\, . 
\end{equation}
Notice that \eqref{cl} and $n_0 \in L^1_+ $ guarantees nonnegative solutions $n(t,x)\ge 0$ and mass conservation (at least formally):
\begin{equation}\label{masscons}
\int_{x>0} n(t,x)\,\D x =\int_{x>0} n_0(x)\,\D x = M\, . 
\end{equation}

This model has been proposed in \cite{Voituriez} to describe basically the polarisation of cells.
The interesting feature of (\ref{eq}) is that the solution may become unbounded in finite time. Such a behavior is called \textit{blow-up in finite time}. In this note we shall prove the following simple alternative:\medskip

\begin{theorem}\label{thm:critical mass} Assume $n_0(x)$ is continuous on $[0,+\infty)$ and $n_0\in L^1_+((1+x)dx)$. 
If $M\leq1$ the solution of \eqref{eq}--\eqref{cl} is global in time. On the contrary if $M>1$, assume in addition that $n_0$ is non increasing, then the solution of  \eqref{eq}--\eqref{cl} blows-up in finite time.
\end{theorem}\medskip

\begin{remark} It would be possible to weaken the assumptions on $n_0(x)$ (basically $\int_{x>0} n_0(x)\left|\log n_0(x)\right|\, \D x<+\infty$) by using strong regularizing effects of the laplacian (at least in the subcritical case $M<1$) but this is beyond the scope of this note.
\end{remark}\medskip

\begin{remark}
Such a critical mass phenomenon (global existence {\em versus} blow-up depending on the initial mass) has been widely studied for the Keller-Segel (KS) system: $\partial _t n(t,x) = \partial _{xx} n(t,x) - \partial _x( n(t,x)\partial_x c(t,x))$ (also known as the Smoluchowski-Poisson system) in two dimensions of space (see \cite{BDP} and references therein). The KS system describes macroscopically a population of diffusive particles which attract each other through a diffusive chemical signal (resp. gravitational field), solution of the Poisson equation: $-\Delta c(t,x) = n(t,x)$ with homogeneous Neumann boundary conditions \cite{Keller.Segel,JL,Nagai,HorstmannI}. On the other hand the chemical field in \eqref{eq} is in fact solution of the Laplace equation with inhomogeneous Neumann boundary conditions: $-\partial_x c(t,0) = n(t,0)$ (production of the signal is located on the boundary). Although the Keller-Segel cannot exhibit blowing-up solutions in one dimension of space, it is indeed the case for \eqref{eq} (Theorem \ref{thm:critical mass}). As a conclusion, \eqref{eq} appears to have the same "singularity" as   the two-dimensional Keller-Segel system. 
Note that there exist other ways to mimick the two dimensional case singular behaviour of KS in one dimension \cite{BilerWoyczynski,CPS,CieslakLaurencot}. 
\end{remark}\medskip

\begin{remark}
There is a strong connection between the equation under interest here \eqref{eq} and the one-dimensional Stefan problem. The later writes indeed \cite{HV1}:
\begin{equation}
\left\{\begin{array}{l}
\partial_t u(t,z) = \partial_{zz} u(t,z) \, , \quad \, t>0\, , \, z\in (-\infty,s(t))\, , \\
\lim_{z\to -\infty}\partial_zu (t,z) = 0 \, , \quad u(t,s(t)) = 0\, , \quad \partial_z u (t,s(t)) = -s'(t)\, .
\end{array}\right.
\end{equation}
The temperature is initially nonnegative: $u(0,z) = u_0(z)\geq 0$. By performing the following change of variables: $\phi(t,x) = - u(t,s(t)-x)$, we get the following equation:
\begin{equation}
\left\{\begin{array}{l}
\partial_t \phi(t,x) = \partial_{xx} \phi(t,x) - s'(t) \partial_x\phi(t,x)\, , \quad \, t>0\, , \, x\in (0,+\infty)\, , \\
\lim_{x\to +\infty}\partial_x\phi (t,x) = 0 \, , \quad \phi(t,0) = 0\, , \quad \partial_x \phi (t,0) = -s'(t)\, .
\end{array}\right.
\end{equation}
By differentiating this equation, we recover \eqref{eq} for $n(t,x) = \partial_x \phi(t,x)$. The condition $\phi(t,0) = 0$ turns out to be the mass conservation of $n(t,x)$.

This connection provides some insights concerning the possible continuation of solutions after blow-up \cite{HV1}. This question has raised a lot of interest in the past recent years \cite{HV2,V1,V2,DS}. It is postulated in \cite{HV1} that the one-dimensional Stefan problem is generically non continuable after the blow-up time.

\end{remark}\medskip

Using ad-hoc entropy methods (which are to be adapted to the nonlinearity in this problem), we are able to investigate long-time behaviour in the critical ($M=1$) and the subcritical case ($M<1$): this is the purpose of Theorems \ref{th:long time critical} and \ref{th:long time subcritical}.
In short, the results read as follows: there exists a one-parameter family of stationary states for the critical mass only  (namely decreasing exponentials). In this case the conservation of the first momentum enables to select one particular profile among this family. 
In the subcritical case, an appropriate rescaling has to be performed in order to capture the intermediate asymptotics. For each mass $M<1$ there exists a unique stationary state (with explicit formulation), and we prove convergence (in relative entropy) of the rescaled solution towards this profile (namely the product of a decreasing exponential and a Gaussian profile). The results are clearly similar to the classical Keller-Segel in two dimensions \cite{BDP}, except that the density converges towards a Dirac mass in the critical case \cite{BCM}.

\section{The critical mass phenomenon}

\paragraph*{Blow-up for $M>1$.}
To prove that solutions blow-up in finite time, we show that the first momentum of $n(t,x)$ cannot remain positive for all time. This technique was first used by Nagai \cite{Nagai}, then by many authors in various contexts (see \cite{Biler95,BilerWoyczynski,Corrias.Perthame.Zaag,DP,CieslakLaurencot} for instance). Other strategies have been used to prove the existence of blowing-up solutions (either constructive by Herrero and Velazquez \cite{HV2} or undirect \cite{Horstmann.Wang01}), however up to date this trick is the only way to provide explicit criterion and appears to be quite robust to variations around Keller-Segel \cite{Blanchet.Carrillo.Laurencot,Calvez.Corrias.Ebde}.

First, the assumption that $n_0$ is a nonincreasing function guarantees that $n(t,\cdot)$ is also a nonincreasing function for any time $t>0$ due to the maximum principle (notice that the derivative $v(t,x) = \partial_x n(t,x)$ satisfies a parabolic type equation, is initially nonpositive, and is nonpositive on the boundary due to \eqref{cl}).  
Therefore $-\partial_x n(t,x)/n(t,0)$ is a probability density at any time $t>0$. We deduce from Jensen's inequality the following interpolation estimate:
\begin{align*}
\left(\int_{x>0} x \dfrac{-\partial_x n(t,x)}{n(t,0)}\,\D x\right)^2  & \leq 
\int_{x>0} x^2 \dfrac{-\partial_x n(t,x)}{n(t,0)}\,\D x\, , \\
M^2 & \leq 2 n(t,0)  \int_{x>0} x   n(t,x) \,\D x \, . 
\end{align*}
Secondly introduce the first momentum $\bJ(t) = \int_{x>0} x n(t,x)\,\D x$. We have for $M>1$:
\begin{align}
 \frac{\D \bJ(t) }{\D t} & = n(t,0) - M n(t,0) 
 \leq \dfrac{M^2}{2 \bJ(t)} (1 - M) \, , \label{eq:moment1}\\
 \dfrac{\D \bJ(t)^2 }{\D t} & \leq  {M^2} (1 - M)\nonumber\, .
\end{align}
Therefore blow-up of the solution occurs in finite time if $M>1$.

\paragraph*{Global existence for $M<1$.} Global existence results for Keller-Segel type systems have been initiated by J\"ager and Luckhaus \cite{JL}  in the two dimensional case. It relies on a mixture of Gagliardo-Nirenberg type and interpolation inequalities. The novelty here is to use a trace-type Sobolev inequality (simple in the one-dimensional setting) which is required due to the location of the chemical source on the boundary.

We compute the evolution of the density entropy as following,  
\begin{align*}
  \frac{ \D }{ \D t} \int_{x>0} n(t,x)\log n(t,x) \,\D x  & = \int _{x>0} \partial _t n(t,x) \log n(t,x)  \,\D x \\
 & = - \int _{x>0} \left(\partial_x n (t,x) + n(t,0) n(t,x)\right) \dfrac{\partial_x n(t,x) }{n(t,x)}  \,\D x  
 \\& 
 = -  \int _{x>0} \left( \partial_x \log n(t,x) \right)^2 n(t,x)\, \D x   +    n(t,0) ^2 \, .
\end{align*}
The one-dimensional trace inequality we mentioned above writes as following,
\begin{align}
n(t,0) & = - \int_{x>0}  \partial_x n(t,x)\,\D x 
= -  \int_{x>0} \left(\partial_x \log n(t,x)\right)  n (t,x)\,\D x  \, ,
\nonumber\\
n(t,0)^2 &  \leq   M \int_{x>0} \left( \partial_x \log n(t,x) \right)^2 n(t,x)\,\D x  \, .
\label{eq:trace1Da}
\end{align}
Therefore we deduce that 
\begin{equation}
\frac{ \D }{ \D t} \int_{x>0} n(t,x)\log n(t,x) \,\D x  \leq  (M-1)\int_{x>0} \left( \partial_x \log n(t,x) \right)^2 n(t,x)\,\D x \, ,
\label{eq:entropy bound}
\end{equation}
hence the entropy is nonincreasing when the mass is smaller than 1.
Observe that equality holds in the trace inequality \eqref{eq:trace1Da} if $\log n(t,x)$ is constant w.r.t. $x$: there exists $\alpha (t)>0$ such that $n(t,x) = M \alpha(t) \exp(-\alpha(t) x)$. In fact the boundary condition \eqref{cl} implies $M = 1$, which is the only configuration where a stationary state can exist (see Section \ref{sec:StSt}).

A major step towards a complete existence theory of \eqref{eq} in the subcritical is to ensure that the boundary value $n(t,0)$ makes perfect sense. This is a consequence of an Aubin-Lions type argument \cite{AubinLions}, which is straightforward in this over-simplified context.
We ask for continuity w.r.t. $x$ of the density $n(t,x)$:
\begin{equation*}
\left(n(t,x) - n(t,y)\right)^2 \le  \left( \int_{x}^y n(t,z)\, \D z \right) \left( \int_{ x}^y (\partial_x\log n(t,z))^2 n(t,z)\, \D z \right) \, .
\end{equation*}
The bound \eqref{eq:entropy bound} together with the control of moments guarantee that $\int_{ x>0} (\partial_x\log n(t,x))^2 n(t,x)\, \D x $ is finite almost every time. Therefore $n(t,\cdot)$ is continuous almost every time.

To conclude this Section, let us mention that it is now classical to prove suitable regularizing effects acting on \eqref{eq} in the subcritical case $M<1$. Indeed an {\em a priori} estimate \eqref{eq:entropy bound} on $\int_{x>0} n(t,x)\left(\log n(t,x)\right)_+ \,\D x $  yields the boundedness of all $L^p-$norms ($1<p<+\infty$) \cite{JL,CC,BDP,CPS}.

\section{Long-time behaviour: convergence in relative entropy for the critical and the subcritical cases}
\label{sec:StSt}
\paragraph*{The critical case.} Equilibrium configurations for the cell density are only possible when the mass is critical: $M=1$ (as it is for the two-dimensional Keller-Segel problem). In this case, a straigtforward computation leads to the one-parameter family:
\begin{equation}\label{steadystate}
h_\alpha(x)=\alpha e^{-\alpha x}\, , \quad \alpha>0 .
\end{equation} 
On the other hand, notice that the first momentum of the cell density is conserved \eqref{eq:moment1}. This prescribes a unique choice for $\alpha$: 
$ \alpha^{-1} =  \bJ(0)\, $.

\begin{theorem}\label{th:long time critical}
Assume $n_0(x)$ being as in Theorem \ref{thm:critical mass}, and the mass being critical: $M = 1$.
As time goes to infinity, the cell density converges (in relative entropy) towards $h_\alpha(x)$.
\end{theorem}\medskip

The convergence proof is based on evaluating the time evolution of the relative entropy, defined as follows:
\begin{equation} \label{eq:relative entropy}
\bH(t)= \int_{x>0}\frac{n(t,x)}{h_\alpha (x)} \log \left(\frac{n(t,x)}{h_\alpha (x)}\right)h_\alpha (x)\, \D x\, . 
\end{equation}
The precise description of the equality cases for inequality \eqref{eq:trace1Da} enables to perform accurate estimates.
A direct computation yields the following estimate:
\begin{equation}
\frac{ \D }{ \D t} \bH(t) = - \int_{x>0} \left(\partial_x \log n(t,x) + n(t,0)\right)^2 n(t,x)\, \D x\, .
\end{equation}
We refer to \cite{CalvezMeunierVoituriez} for more details.

\paragraph*{Self-similar decay in the subcritical case.}
In the sub-critical case $M<1$ one expects the density $n(t,x)$ to decay self-similarly. For this purpose the density is appropriately rescaled:
\[ n(t,x) = \dfrac{1}{\sqrt{1+2t}} u\left( \log\sqrt{1+2t}  , \dfrac{x}{\sqrt{1+2t}}\right)\, . \]
The new density $u(\tau,y)$ satisfies:
\begin{equation}
\partial_\tau u(\tau,y) = \partial_{yy} u(\tau,y)   + \partial_y \left(y u(\tau,y) \right) + u(\tau,0)\partial_y u(\tau,y) \, ,
\end{equation}
and no-flux boundary conditions: $\partial _y u (\tau,0)+u(\tau,0)^2 =0$.
The additionnal left-sided drift contributes to confine the mass in the new frame $(\tau,y)$. The stationary equilibrium in this new setting can be computed explicitely. The expected self-similar profile writes:
$g_\alpha(y) = \alpha\exp\left(-\alpha y - {y^2}/2\right)$,  where $\alpha$ is given by  the relation $P(\alpha) = M$, $P$ being an increasing function defined as following:
\[ P(\alpha) = \int_{y>0}  \exp\left(-  y - \dfrac{y^2}{2\alpha^2}\right)\, \D y \, , \quad \left\{\begin{array}{r} \lim_{\alpha \to 0} P(\alpha) = 0 \\ \lim_{\alpha \to +\infty} P(\alpha) = 1 \end{array}\right. \,   . \]

\begin{theorem}\label{th:long time subcritical}
Assume $n_0(x)$ being as in Theorem \ref{thm:critical mass}, and the mass being subcritical: $M < 1$.
As time goes to infinity, the first momentum $\bJ(\tau)$ of the density converges to $\alpha(1-M)$ and the cell density converges (in relative entropy) towards $g_\alpha(y)$.  
\end{theorem}\medskip

The proof of this Theorem relies again on the time evolution of the relative entropy:
\begin{equation} 
\bH(\tau) = \int_{y>0}\frac{u(\tau,y)}{g_\alpha (y)} \log \left(\frac{u(\tau,y)}{g_\alpha (y)}\right)g_\alpha (y)\, \D y\, .  
\end{equation}
More precisely we have:
\begin{multline}
 \frac{\D }{\D \tau} \left\{\bH(\tau)  + \dfrac{1}{2(1-M)}  \left( \bJ(\tau) - \alpha(1-M) \right)^2\right\}  \\ =  - \int_{y>0} u(\tau,y) \left(   \partial_y \log u(\tau,y)  + y + u(\tau,0) \right)^2\, \D y  - \dfrac1{(1-M)}\left(\frac{\D }{\D \tau} \bJ(\tau)\right)^2\, . 
\end{multline} 

\section{Analysis of a coupled ODE/PDE model}

We investigate in this section a variant of \eqref{eq} which is more relevant for modelling purposes \cite{Voituriez}. In this new setting, the chemical is supplied by a quantity $\mu(t)$ which evolves by exchanging particles at the boundary $x=0$: 
\begin{equation} \label{eq:Raphael1D}
\left\{\begin{array}{l} 
\partial _t n (t,x)=\partial _{xx} n (t,x) + \mu (t) \partial _x n (t,x) \, , \quad t >0\, , \, x\in (0,+\infty)\, , \\ 
 \mu'(t) = n(t,0)- \mu(t)\, ,
\end{array}\right.
\end{equation}
together with the initial conditions: $n(t=0,x) = n_0(x)\geq 0$ and $\mu(t=0) = \mu_0$. The conservation of the total mass of particles:
\begin{equation} \int_{x>0} n(t,x)\, dx + \mu(t) = M\, , \label{eq:mass conservation}\end{equation} 
yields the following boundary condition for the cell density:\[ \partial _x n (t,0)+\mu(t) n (t,0)=\mu'(t) \, . \label{eq:BC2}\] 


\paragraph*{Long-time convergence in the case $M>1$.} We denote by $m(t)$ the mass of the cell density $n(t,x)$: 
\begin{equation}
m(t) = \int_{x>0} n(t,x)\, \D x\, .
\end{equation}
(notice $m'(t) + \mu'(t) = 0 $ due to the conservation of mass).
Introduce the relative entropy:
\[
\bH(t)= \int_{x>0}\frac{n(t,x)}{m(t)h  (x)} \log \left(\frac{n(t,x)}{m(t)h(x)}\right)h (x)\, \D x\, ,
\]
where the expected profile $h$ is given by:
\[ h(x) = \overline{\mu} \exp\left(-\overline{\mu} x\right)\, , \quad \overline{\mu} = M-1 \, . \]

\begin{theorem}
As time goes to infinity, the mass $m(t)$ of the cell density converges to 1 and the cell density converges (in relative entropy) towards $h(x)$.
\end{theorem}\medskip

The  proof of this Theorem relies again on the time evolution of the relative entropy. This is strongly inspired from the previous computation, but takes into consideration the non-conservation of mass for the cell density and the dynamics of $\mu(t)$ \cite{CalvezMeunierVoituriez}.

\bibliographystyle{plain}
\def\cprime{$'$} \def\lfhook#1{\setbox0=\hbox{#1}{\ooalign{\hidewidth
  \lower1.5ex\hbox{'}\hidewidth\crcr\unhbox0}}} \def\cprime{$'$}
  \def\cprime{$'$} \def\cprime{$'$} \def\cprime{$'$} \def\cprime{$'$}

\noindent{\em Acknowledgement: J.~Van Schaftingen, J.J.L.~Vel\'azquez, R.~Voituriez, CRM.}

\end{document}